\documentclass[11pt,reqno]{article}
\usepackage{amsmath, latexsym, amsfonts, amssymb,
amsthm, amscd,epsfig,enumerate}
\advance\hoffset -.75cm

\usepackage{amsfonts}
\usepackage{latexsym}
\usepackage{amsmath}
\usepackage{amssymb}

\newcommand{\C}{\mathbb C}

\newcommand{\Z}{\mathbb Z}
\newcommand{\N}{\mathbb N}
\newcommand{\Prob}{\mathbb P}

\newcommand{\eps}{\varepsilon}

\newcommand{\cE}{\mathcal{E}}

\newcommand{\cF}{\mathcal{F}}

\newcommand{\pr}{\mathbb P}
\newcommand{\T}{\mathbb T}

\newcommand{\ident}{{\mathchoice {\rm 1\mskip-4mu l} {\rm 1\mskip-4mu l}
{\rm 1\mskip-4.5mu l} {\rm 1\mskip-5mu l}}}

\newtheorem{teo}{Theorem}[section]
\newtheorem{lem}[teo]{Lemma}
\newtheorem{cor}[teo]{Corollary}
\newtheorem{rem}[teo]{Remark}
\newtheorem{pro}[teo]{Proposition}
\newtheorem{defn}[teo]{Definition}

\newtheorem{teo2}{Theorem}
\newtheorem{exm}[teo2]{Example}

\title
{Weak survival for branching random walks on graphs}

\author
{Daniela Bertacchi\\
Dipartimento di Matematica e Applicazioni\\
Universit\`a di Milano--Bicocca\\
via Cozzi 53, 20125 Milano, Italy\\
daniela.bertacchi\@@unimib.it
\and
Fabio Zucca \\
Dipartimento di Matematica \\
Politecnico di Milano\\
piazza Leonardo da Vinci 32, 20133 Milano, Italy\\
fabio.zucca\@@polimi.it}

\date{}

\begin{document}

\maketitle

\begin{abstract}
We study weak and strong survival for branching random walks on
multigraphs.
We prove that, at the strong critical value, the process dies out locally
almost surely. We relate the weak critical value to
a geometrical parameter of the multigraph.
For a large class of multigraphs we prove that,
at the weak critical value, the process dies out globally
almost surely. Moreover for the same class we prove that
the existence of a pure weak phase is equivalent to nonamenability;
this improves a result of Stacey \cite{cf:Stacey03}.
\end{abstract}



\baselineskip .6 cm

\section{Introduction}\label{sec:intro}
\setcounter{equation}{0}

In recent years, much study has been devoted to various stochastic processes,
such as percolation, Ising model, contact process and
branching random walk, on general graphs
(\cite{cf:Hagg1}, \cite{cf:Lyons1}, \cite{cf:Lyons2}, \cite{cf:Pem},
\cite{cf:PemStac1} only to mention a few,
see \cite{cf:Lyons3} for more references).
A double motivation underlies the search for settings
other than the usual $\Z^d$: on one hand the need for structures which may
serve as models for inhomogeneous crystals, biological
structures or social networks and on the
other hand the fact that on general graphs interesting phenomena,
which are absent in $\Z^d$, are observed.
In particular the branching random walk (BRW) has been
studied on trees (see \cite{cf:MadrasSchi}, \cite{cf:Ligg1}, \cite{cf:Ligg2},
\cite{cf:HuLalley}, \cite{cf:PemStac1}) and on quasi-transitive graphs
(see \cite{cf:Stacey03}).

In this paper we study the BRW  on a connected
multigraph $X$ with bounded degree (see Section \ref{subse:def} for the formal definition).
Roughly speaking a $\lambda$-BRW can be described by the following rules: each particle dies after
an exponential time with parameter 1 
and breeds independently on each neighbor at exponential intervals
with parameter $\lambda$.
We start with a finite number of particles, hence the $\lambda$-BRW can be viewed as
a continuous-time random walk on the countable state space of finite configurations
$\eta \in \N^X$. On each site $x \in X$ the transitions are:
\[
\begin{split}
\eta(x) \to \eta(x)-1 &\qquad \hbox{at rate } \eta(x), 
\\
\eta(x) \to \eta(x)+1 &\qquad \hbox{at rate } \lambda \sum_{y \in D(x)} \eta(y), \\
\end{split}
\]
where $\eta(x)$ is the number of particles at site $x$ and $D(x)$ is the set of neighbors of $x$
(see Section \ref{subse:def}).

The BRW has originally been introduced as a model for biological populations
dynamics (although it has been argued that this model is far from being
satisfactory, see for instance the discussion in \cite{cf:Jagers})
and, besides being interesting in itself, has also been studied
for its relationship with the contact process (the process which has the same transition rules of the BRW but
state space $\{0,1\}^X$). Indeed the BRW stochastically dominates the contact process and has an additive property which the
contact process lacks: the sum of two $\lambda$-BRWs is still a $\lambda$-BRW.

The $\lambda$-BRW on $\Z^d$ shows only two possible behaviors (called \textit{phases}):
if $\lambda  \le 1/2d$ there is extinction almost surely; if $\lambda  > 1/2d$, for all $t_0>0$ we have that
$\pr(\eta_t(0)>0  \hbox{ for some } t \ge t_0 )>0$ where $\eta_t(0)$ is the number of particles at $0$ at time $t$.
The main interest of the study of
BRW on trees is that a third phase appears. Indeed
 we may
identify two kinds of survival:
\begin{enumerate}[$(i)$]
\item
\textit{weak} (or global) \textit{survival} -- the total number of particles is positive at each time;
\item
\textit{strong} (or local) \textit{survival} -- the number of particles at one site $x$ is not eventually $0$.
\end{enumerate}
In the first case it is easy to see that the total number of particles diverges (see Section \ref{sec:main} for details);
in the second case the survival at a site $x$ does not depend on the site chosen.

Let us denote by $\lambda_w$ (resp.~$\lambda_s$) the infimum of the values $\lambda$ such that there is weak (resp.~strong)
survival. Clearly $\lambda_w \le \lambda_s$ and we may have three distinct phases corresponding to the following intervals for $\lambda$:
$[0, \lambda_w)$, $(\lambda_w,\lambda_s)$, $(\lambda_s,+\infty)$. The middle interval may be empty;
if, on the contrary, $\lambda_w < \lambda_s$ then we say that
the BRW has a pure weak phase. In this phase the process leaves any finite subset eventually a.s., hence it
survives globally by drifting to infinity (see \cite{cf:HuLalley} for details on the convergence to the boundary in the case of
homogeneous trees).

This paper is devoted to three main issues: the identification of the critical value $\lambda_w$, the behavior of the process at the critical values $\lambda=\lambda_s$
and $\lambda=\lambda_w$
and the existence of the pure weak phase.
In \cite{cf:PemStac1} it was proved that $\lambda_s$ is related to a particular asymptotic
degree of the graph. We prove that, under some
rather general geometrical conditions on the multigraph, $\lambda_w$ is related to
another asymptotic degree (Theorems \ref{th:conditioncones}, \ref{th:Tnk} and \ref{th:almostregular}).
Moreover we prove, by using generating functions techniques, that if $\lambda=\lambda_s$ the
process dies out locally a.s.~(Theorem \ref{th:pemantleimproved})
and that, if $\lambda=\lambda_w$, on  a large class of multigraphs
the process dies out globally a.s.~(Theorem \ref{th:almostregular}).
The use of multigraphs is mainly needed in view of Definition~\ref{Fclass} that defines the
class of (multi)graphs for which our results hold.

As for conditions for the existence of the pure weak phase, one is lead to investigate 
nonamenable graphs. Indeed, usually, nonamenable graphs are graphs where certain phenomena,
absent in the amenable case, appear (see \cite{cf:Lyons3} for a survey). Nevertheless a statement
like ``nonamenability of the graph is equivalent to the existence of a pure weak phase for
the BRW'' has been disproved in \cite{cf:PemStac1}. The authors showed a nonamenable tree
where the BRW has no pure weak phase and an amenable tree where there is such a phase (note that
these counterexamples are both of bounded degree). Hence one hopes to prove a similar statement
for a more restricted class of graphs. Work is this direction has been done in
\cite[Theorem 3.1]{cf:Stacey03} which states the equivalence between nonamenability and the
existence of a pure weak phase for quasi-transitive graphs. We prove the same equivalence
for a larger class of multigraphs which strictly includes both quasi-transitive graphs and regular
graphs (Theorem \ref{th:nonam} and Example \ref{exmp0}).

Let us give the outline of the paper. In Section \ref{sec:def} we introduce
the main definitions and 
we define some
generating functions and a generalized branching process which will be useful in the sequel.
Moreover we introduce two \textit{asymptotic degrees}
$M_s$ and $M_w$ which depend only on the geometrical structure of the multigraph.

Section \ref{sec:main} is devoted to the detailed study of $\lambda_s$ and $\lambda_w$.
We give a sufficient condition for the absence of the pure weak phase,
which, in particular,  implies that there is no weak phase on polynomially growing multigraphs
(Corollary \ref{cor:1}).
 We recall the well-known characterization of
$\lambda_s=1/M_s$ and we show that for a large class of multigraphs $\lambda_w=1/M_w$ (see Theorem \ref{th:conditioncones}).
Clearly for this class we have that $\lambda_w<\lambda_s$ if and only if $M_s < M_w$.
We give two different sufficient conditions for a multigraph to satisfy the hypotheses of Theorem \ref{th:conditioncones}.
The first condition (Theorem \ref{th:Tnk}) is satisfied, for instance, by certain radial trees which are not
quasi transitive; for these trees we show that nonamenability is equivalent to the existence of the pure weak phase
(see Example \ref{exmp0}). As for the second condition (Theorem \ref{th:almostregular}), we introduce a class of morphisms
(see Definition \ref{Fclass}) of multigraphs and we show that it preserves $\lambda_w$, $M_w$ and, in some cases,
$\lambda_s$ (Proposition \ref{th:twoprocesses}). By using these morphisms, the class of $\cF$-multigraphs
is defined; for this class we show that the $\lambda_w$-BRW dies out globally almost surely.
Finally, Theorems \ref{th:almostregular} and \ref{th:nonam} yield, for non-oriented
$\cF$-multigraphs, the equivalence of the following conditions:
(i) $\lambda_w < \lambda_s$, (ii) $M_s<M_w$ and (iii) nonamenability.
In Section \ref{sec:examples}
some examples of multigraphs, which can be studied via our results, are given.

The BRW studied in Section \ref{sec:main} may be viewed as a population which reproduces
following an ``edge breeding'' pattern, while some authors prefer a ``site breeding'' pattern.
In Section \ref{sec:modif} we consider this modification of the BRW.
These two versions of the BRW are essentially equivalent on regular graphs, while in
the general setting the behavior of the ``site breeding'' one can be much more easily characterized (see Theorem \ref{th:modif}).
We show that BRWs and modified BRWs may both be seen as particular cases
of BRWs on weighted graphs. Most of
the results given in the previous sections still hold in this general setting.

Section \ref{sec:fr} is devoted to a final discussion of 
open questions.

\section{Basic definitions and preliminaries}\label{sec:def}

\subsection{Multigraphs}\label{subse:def}

A countable (or finite) multigraph is a couple $(X,E(X))$, where
$X$ is the countable (or finite) set of vertices and $E(X) \subseteq X
\times X \times \N_*$ is the set of (oriented) edges (where $\N_*$ is the set of positive natural numbers);
we define the number of edges from $x$ to $y$ as
$n_{xy}:=|\{i:(x,y,i) \in E(X)\}| \equiv \max \{i:(x,y,i) \in E(X)\}$ (where
$|\cdot|$ denotes cardinality). We denote by
$D(x):=\{y \in X: n_{xy}>0\}$ the set of neighbors of $x$ and by
$\deg(x):=\sum_{y \in D(x)} n_{xy}$ the degree of $x$.
 If $n_{xy}=n_{yx}$ for all $x,y \in X$ then the multigraph is
called \textit{non oriented}. A multigraph is a graph if and only if
$n_{xy}=\ident_{D(x)}(y)$.

A path from $x$ to $y$ of length $n$ is a couple of sequences
$\left (\{x=x_0, x_1, \ldots, x_n=y\}, \{k_1, \ldots, k_n\} \right )$ such that
$n_{x_i x_{i+1}} \ge k_{i+1}>0$ for all $i=0,1,\ldots,n-1$. The multigraph is said to be
\textit{connected} if there exists a path (of suitable length) from $x$ to $y$, for all $x,y \in X$.
From now on, the
multigraph will always be connected and of \textit{bounded degree},
that is, $M(X):=\sup_{x \in X} {\rm deg}(x) <+\infty$;
obviously $M$ depends on $(X,E(X))$, nevertheless
to avoid cumbersome notation the dependence on the
set of edges will be tacitly understood. The same implicit assumption
will be made for all quantities depending on the multigraph.
Moreover if not explicitly stated, the multigraph does not need to be non oriented.

Let $\gamma_{x,y}^n$
be the number of paths of length $n$ from $x$ to $y$ (and
$\gamma^0_{x,y}:=\delta_{x,y}$). More explicitly to each
sequence $\{x=x_0, x_1, \ldots, x_n=y\}$ there corresponds a set of
$\prod_{i=0}^{n-1} n_{x_i x_{i+1}}$ paths in the multigraph,
whence $\gamma^n_{x,y}$ is the sum over all the sequences
$\{x=x_0, x_1, \ldots, x_n=y\}$ of $\prod_{i=0}^{n-1} n_{x_i x_{i+1}}$.
Moreover let $T_x^n$ be the number of paths from $x$ of length
$n$, that is, $T_x^n:=\sum_{y \in X} \gamma^n_{x,y}$. Finally, let
$\phi^n_{x,y}$ the number of paths of length $n$ starting from $x$
and reaching $y$ for the first time;
to be precise, $\phi^n_{x,y}$ is the number of paths
$\left (\{x=x_0, x_1, \ldots, x_n=y\}, \{k_1, \ldots, k_n\} \right )$ such that
$x_i \not = y$ for all $i=1,\ldots, n-1$. By definition
$\phi^0_{x,y}:=0$ for all $x,y \in X$.

For $\gamma^n_{x,y}$ and $T_x^n$ the following recursive
relations hold for all $n,m \geq 0$
\[
\begin{cases}
\gamma^{n+m}_{x,y}=\sum_{w \in X} \gamma^n_{x,w} \gamma^m_{w,y} \\
\\
\gamma^1_{x,y}=n_{xy}
%
\end{cases}
\qquad
\begin{cases}
T_x^{n+m}=\sum_{w \in X} \gamma^m_{x,w} T_w^n \\
\\
T_x^1={\rm deg}(x)\\
\end{cases}
\]
and, for all $n \ge 1$,
\[
\gamma_{x,y}^n=\sum_{i=0}^n \phi_{x,y}^i \gamma^{n-i}_{y,y}.
\]

Given any vertex $x \in X$ and $n \in \N$, we define $\rho(x,y):=\min
\{i: \gamma^i_{x,y}>0\}$ and $B(x,n):=\{y
\in X: \rho(x,y) \le n\}$;
note that $\rho$ is a metric if $n_{xy}>0$ is equivalent to $n_{yx}>0$ for
all $x,y \in X$ (for instance in the case of non-oriented multigraphs).

By using the number of paths it is possible to introduce two
asymptotic degrees, namely
\[
M_s(X):=\limsup_{n} (\gamma^n_{x,y})^{1/n} \qquad
M_w(X):=\limsup_{n} (T_x^n)^{1/n}.
\]
It is easy to show that the above definitions do not depend on
the choice of $x,y \in X$, moreover simple arguments of
supermultiplicativity show that $M_s(X)=\lim_{n \to \infty}
(\gamma^{dn}_{x,x})^{1/dn}=\sup_{n} (\gamma^{dn}_{x,x})^{1/dn}$
where
$d:={\rm gcd} \{n: \gamma_{x,x}^n >0 \}$ is the period of the multigraph (which does not depend on
the choice of $x$). Analogously $M_s(X)=\lim_{n \to \infty}
(\gamma^{dn+i}_{x,y})^{1/(dn+i)}$, where $0 \le i \leq d-1$ is uniquely
chosen such that $\gamma^{n}_{x,y}>0$ implies $n=i$ (mod $d$).
In the rest of the paper, whenever there is no ambiguity, we will denote
$M(X)$, $M_s(X)$ and $M_w(X)$ simply by $M$, $M_s$ and $M_w$.

By definition $1 \le M_s \leq M_w \leq M$. We note that $M_w=M$ if the
multigraph is regular, that is, it has constant degree. Moreover
if $|B(x,n)|^{1/n} \to 1$ when $n \to +\infty$ then $M_w=M_s$ (see
Corollary \ref{cor:1}).
It is well known that, for a regular non-oriented graph, $M_s<M_w$ if and only if
it is nonamenable (see Section \ref{sec:nonamenability} for the definition).

\subsection{Generating functions}

In order to find some characterizations of $M_s$ and $M_w$,
let us define the
generating functions
\[
\begin{split}
H(x,y|\lambda)&:=\sum_{n =1}^\infty \gamma^n_{x,y} \lambda^n,
\qquad \Theta(x|\lambda):=\sum_{n =1}^\infty T_x^n \lambda^n,
\end{split}
\]
with radius of convergence $1/M_s$ and $1/M_w$ respectively. Of
course for all $\lambda \in \C$ such that $|\lambda|<1/M_w$ we have
$\Theta(x|\lambda)=\sum_{y \in Y} H(x,y|\lambda)$ and the
following relations hold
\begin{equation}\label{eq:HTheta}
\begin{split}
H(x,y|\lambda)&=\delta_{x,y}+\lambda \sum_{w \in X} \gamma_{x,w}^1
H(w,y|\lambda)\\
& =
\delta_{x,y}+\lambda \sum_{w \in X} H(x,w|\lambda)\gamma_{w,y}^1,
\qquad \forall \lambda \in \C: |\lambda|< 1/M_s,\\
\Theta(x|\lambda)&=1+\lambda \sum_{w\in X} \gamma_{x,w}^1 \Theta(w|\lambda),
\qquad \qquad \forall \lambda \in \C: |\lambda|< 1/M_w.\\
\end{split}
\end{equation}
We define
\[
\Phi(x,y|\lambda):=\sum_{n =1}^\infty \phi_{x,y}^n \lambda^n;
\]
it is easy to see that $\Phi(x,x|\lambda)= \lambda \sum_{y \in X, y \not = x} \gamma^1_{x,y} \Phi(y,x|\lambda)+
\lambda \gamma^1_{x,x}$ and
if $x,y,w \in X$ are distinct vertices such that every path from $x$ to $y$ contains $w$ then
$\Phi(x,y|\lambda)=\Phi(x,w|\lambda)\Phi(w,y|\lambda)$.
Moreover
\[
H(x,y|\lambda)=\Phi(x,y|\lambda)H(y,y|\lambda)+\delta_{x,y}, \quad \forall
\lambda: |\lambda|< 1/M_s.
\]
Since the radius of the series $H(x,x|\cdot)$ does not depend on
the choice of $x \in X$ and since
\begin{equation}\label{eq:genfun1}
H(x,x|\lambda)=\frac{1}{1-\Phi(x,x|\lambda)},
\qquad \forall \lambda \in \C: |\lambda|< 1/M_s,
\end{equation}
we have that $1/M_s=\max\{ \lambda \geq 0 :\Phi(x,x|\lambda)\leq 1\}$
for all $x \in X$ (remember that $\Phi(x,x|\cdot)$ is left-continuous on $[0,1/M_s]$
and that $1/(1-\Phi(x,x|\lambda))$ has no analytic prolongation in $1/M_s$).

The computation of $M_w$ is not easy in general, but in the
case of finite multigraphs
there is a simple characterization of  $M_w$. In the following theorem, $\rm{\mathbf{Id}}$
is the identity matrix.

\begin{teo}\label{th:fromPFTheorem}
Let $(X,E(X))$ be an irreducible, finite multigraph with adjacency matrix
$N:=(n_{xy})_{x,y \in X}$, then
\[
1/M_w=1/M_s=\min\{\lambda>0: {\rm det}(\lambda N -\rm{\mathbf{Id}})=0\}.
\]
\end{teo}
\begin{proof}
We use the same notation $N$ for the matrix and the
linear operator. By the Perron-Frobenius Theorem there exists an
eigenvalue $\sigma_0>0$ of $N$ such that any other eigenvalue
$\sigma$ satisfies $|\sigma|<\sigma_0$ and the same holds for $N^t$. Moreover ${\rm
dim}({\rm Ker}(N^t-\sigma_0\rm{\bf{Id}}))=1$ and it is possible to choose
the eigenvector $v$ in such a way that $v>0$. It is clear that any vector
$w<0$
cannot possibly belong to ${\rm
Rg}(N-\sigma_0{\rm{\bf{Id}}})\equiv {\rm{Ker}}( N^t-\sigma_0 \bf{\rm{Id}})^\bot$
since
$\langle w,v\rangle <0$. Then the equation \eqref{eq:HTheta} (which holds for $|\lambda|<1/M_w$),
can be written as
\begin{equation}\label{eq:system}
(\lambda N-\rm{\bf{Id}})\Theta(\lambda)=-\mathbf{1}=:
\begin{pmatrix}
1 \\
1 \\
\vdots \\
1 \\
\end{pmatrix},
\end{equation}
and has no solutions if $\lambda=1/\sigma_0$. On the other hand
equation \eqref{eq:system} defines a holomorphic (vector) function 
$\Theta^\prime(\lambda)=(\lambda N-\rm{\bf{Id}})^{-1}\mathbf{1}$ on
$\{\lambda \in \C: |\lambda|< 1/\sigma_0\}$. Note that $\Theta^\prime$ coincides with $\Theta$
on $\{\lambda \in \C: |\lambda|< \min\{1/\sigma_0,1/M_w\}\}$, hence 
$1/\sigma_0\le 1/M_w$. If $1/\sigma_0<1/M_w$ then 
there would be an analytic prolongation of $\Theta^\prime$ to $1/\sigma_0$ and
by continuity eq.~\eqref{eq:system} would hold for $\lambda=1/\sigma_0$.
\end{proof}

\subsection{Generalized branching process}\label{sec:gbp}

In the classical branching process (see for instance \cite{Harris}) there is a unique offspring
distribution according to which each individual breeds. We consider a
generalized branching process where each father may have different
types of children and each of them
breeds according to a
specific distribution which depends on its type and on the father.
To be more specific, let $\mathbb
T=(\bigcup_{i=0}^\infty\N_*^{2i}, E(\mathbb T))$ where $\N_*^0:=\{o\}$ where
$o$ is the root of the tree $\mathbb T$. 
Identifying as usual $\N_*^{2n}\times \N_*^2$ with $\N_*^{2n+2}$,
the set of edges is
\[
E(\mathbb T):=\{(x,y)\in\mathbb T:\exists k\in \N_*^2,
y=(x,k)\}\cup \{(o,k):k\in\N_*^2\}.
\]
Roughly speaking $y=(x,i,j)$ means that $y$ is the $j$-th son of
type $i$ of its father $x$ (whereas $(i,j)$ is the $j$-th son of
type $i$ of $o$) and the oriented edges are drawn from fathers to sons. 
Moreover $\bigcup_{i=0}^n \N_*^{2i}$ represents the genealogic tree
of the progenies of $o$ up to the $n$-th generation.
We provide each individual $x$ with a distribution $\mu_x$ such that
if $x=(v,i,j)$ and $y=(v,i,k)$ then $\mu_x\equiv\mu_y$ (that is,
the offspring distribution depends only on the father and on the type).
Now, each distribution is defined on the countable space
$\cE:=\{f\in\N^{\N_*}: S(f)<+\infty\}$ where
$S(f)=\sum_{i=1}^\infty f(i)$.
To be more precise
it is possible to construct a canonical probability space $(\Omega,\mathcal A,\pr)$
supporting the generalized branching process and such that $\pr$ satisfies
\[
\mu_x(f)=\pr \left (
\bigcap_{i=1}^\infty \{x \hbox{ has } f(i) \hbox{ sons of type } i
\}
\right ),
\qquad \forall f \in \cE.
\]
Moreover for every $x\in\mathbb T$ let
$\nu_x$
be the distribution of the total number of children of $x$, that is, 
$\nu_x(k)=\mu_x(\{f:S(f)=k\})$ for
all $k\in\N$.
Take a family of independent $\cE$-valued random
variables $\{Z_x\}_{x\in\mathbb T}$ such that $Z_x$ has
distribution $\mu_x$.

Let us recursively construct this generalized branching process $\{B_n\}_{n\ge0}$:
\[
B_0=\{o\}, \qquad B_{n+1}=\{(v,i,j):v\in B_n, 1\le j\le Z_v(i)\},
\]
where $B_n$ is the $n$-th generation, and its member $v$ has
exactly $Z_v(i)$ children of type $i$. Extinction is
($B_n=\emptyset$ eventually).

\begin{lem}\label{th:powerhouse}
Let $G_x(z)$ be the generating function of $\nu_x$ and suppose
that there exists $\delta\in [0,1)$ such that
$G_x(\delta)\le\delta$ for all $x\in\mathbb T$. Then
$\Prob(B_n=\emptyset\hbox{ {\rm eventually}})\le\delta$.
\end{lem}
\begin{proof}
Denote by $A_n^x$ the event of extinction before the $n$-th generation
of the progenies of $x$.
Let $q_n^x:=\pr(A_n^x)$, clearly $q^x_n$ depends only on the father and the type of $x$;
we claim that $q_n^x\le\delta$ for all $x \in \mathbb T$. We proceed by induction on $n$.
Obviously, for each $x\in \mathbb T$,
$q_0^x=\nu_x(0)
=G_x(0)\le\delta$.
By induction, using the hypothesis of independence,
\[
\begin{split}
q_{n+1}^x & =
\pr \left ( \bigcup_{i=0}^\infty \,\bigcup_{f:S(f)=i} \, \bigcap_{j=1}^\infty \bigcap_{k=1}^{f(j)} A_n^{(x,j,k)}
\right )= \sum_{i=0}^\infty \sum_{f:S(f)=i}
\mu_x(f)\prod_{j=1}^\infty \left (q_n^{(x,j,1)} \right )^{f(j)}\\
& \le
\sum_{i=0}^\infty \sum_{f:S(f)=i} \mu_x(f)\delta^{S(f)}
 = \sum_{i=0}^\infty \delta^i \sum_{f:S(f)=i} \mu_x(f)=
\sum_{i=0}^\infty \delta^i\nu_x(i)=G_x(\delta)\le\delta.
\end{split}
\]
Now, $q_n^o\uparrow\Prob(B_n=\emptyset\hbox{ eventually})$ and
$\delta\ge\lim_n q_n^o$ and this yields the conclusion.
\end{proof}

This lemma trivially applies when each distribution $\mu_x$ is
drawn from a finite set of distributions such that the
corresponding $\nu_x$ represents a supercritical branching
process. In this case we have a finite number of fixed points in $[0,1)$
for the generating functions and
$\delta$ may be taken as the maximum among them
(indeed this is what we do in Theorem \ref{th:conditioncones}).

\section{Main results}\label{sec:main}

\subsection{The critical values}\label{sec:critical}

We investigate the critical values $\lambda_s$ and
$\lambda_w$, their relationship with $M_s$ and $M_w$ and the
behavior of the $\lambda$-BRW when $\lambda=\lambda_s$ or
$\lambda=\lambda_w$.
Since the critical values  do not depend on the number of particles at $t=0$
(nor on their location), we suppose that the initial state is
one particle at a fixed vertex $o \in X$.
To each particle $p$ (present at some time at a site $x$) there corresponds a (unique)
\textit{reproduction trail} starting from the initial particle located
at $o$ at time $0$ reconstructing the genealogy of $p$. 
Roughly speaking, the (space-time) reproduction trail corresponding to $p$
is a path $(\{x_0=o,x_1,\ldots,x_{n-1},x_n=x\},\{k_1,\ldots,k_n\})$ 
along with a sequence $(t_0,\ldots,t_{n-1})$
where $t_0$ is the epoch when the original particle in $o$ generated the ancestor
of $p$ in $x_1$ (through the edge $(o,x_1,k_1)$) and, for $i=1,\ldots,n-1$, 
$t_i$ is the epoch when the ancestor in $x_i$ generated the one in $x_{i+1}$
(through the edge $(x_i,x_{i+1},k_{i+1})$). Clearly, putting $t_{-1}=0$,
for all $i=0,\ldots,n-1$, $t_i-t_{i-1}$ is the realization of an exponential
random variable with rate $\lambda$ (it is tacitly understood that
each ancestor is alive when breeding). Such a trail is said to have length $n$.
For a detailed construction we refer the reader to \cite[Section 3]{cf:PemStac1}
(where what we call reproduction trail is an infection trail).

In \cite[Lemma 3.1]{cf:PemStac1} it was proved that $\lambda_s=1/M_s$ for any
graph. We use a different approach to extend this result to multigraphs;
this approach allows us to study the critical behavior when $\lambda=\lambda_s$.

\begin{teo}\label{th:pemantleimproved}
For each multigraph $(X,E(X))$ we have that $\lambda_s=1/M_s$ and if
$\lambda=\lambda_s$ then the $\lambda$-BRW dies out locally almost surely.
\end{teo}

\begin{proof}
Let us consider a path $\Pi:=\left (\{o=x_0, x_1, \ldots, x_n=o\}, \{k_1, \ldots, k_n\} \right )$ and let
us define its number of cycles $\mathbb{L}(\Pi):=|\{i=1,
\ldots,n:x_i=o\}|$; the expected number of trails along such a path
is $\lambda^n$
(hence to each sequence $\{x_0, x_1, \ldots, x_n\}$ there corresponds a number
$\lambda^n \prod_{i=0}^{n-1} n_{x_i x_{i+1}}$ of expected trails).
Disregarding the original time scale, to the BRW there
corresponds a Galton-Watson branching process: given any particle $p$ in $o$ 
(corresponding to a trail with $n$ cycles), define its children
as all the particles whose trail is a prolongation of the trail of $p$ and is associated 
with a spatial path with $n+1$ cycles.
Hence a particle is of the $k$-th generation if and only if the
corresponding trail has $k$ cycles; moreover it  has one (and only one)
parent in the $(k-1)$-th generation. Since each particle behaves
independently of the others then the process is markovian. Thus the BRW
survives if and only if this branching process does.
The expected number of children of the branching process is the sum
over $n$ of the expected number of trails of length
$n$ and one cycle, that is $\sum_{n=1}^\infty \phi_{o,o}^n\lambda^n=\Phi(o,o|\lambda)$.
Thus we have a.s.~local extinction if and only if $\Phi(o,o|\lambda)\leq 1$, that is,
$\lambda \leq 1/M_s$ (see eq.~\eqref{eq:genfun1} and the remark thereafter).
\end{proof}

Considering the equivalence between a $\lambda$-BRW and a
branching process as discussed in the previous proof, it is clear
that if $\lambda> \lambda_s$, then the conditional probability of
local explosion given non-extinction  is $1$. The same holds
(globally) if $\lambda>\lambda_w$. Indeed the BRW (starting with a
finite number of particles) is a continuous-time random walk on
the countable state space of finite configurations $\eta \in
\N^X$, with a trap state in $\underline 0$ (the configuration with
no particles). Hence all the states but $\underline 0$ are
transient and
 the process which does not hit $\underline 0$ leaves 
$A_k=\{\eta \in \N^X: \sum_{x \in X} \eta(x) \le k\}$
eventually for all $k \in \N$. Indeed, the probability of reaching $\underline 0$ 
starting from any configuration in $A_k$
is uniformly different from $0$ (remember that the reproduction rate is bounded from above in
a bounded degree multigraph), hence the claim follows. 

%
%
%
%
%
%

\noindent Now we focus our attention on the weak critical value.

\begin{lem}\label{lem:general}
For every multigraph
we have that $\lambda_w \geq 1/M_w$.
\end{lem}
\begin{proof}
Since the average number of trails on a fixed path of length
$n$ starting from $(o,0)$  is $\lambda^n$,
the average number of all the trails on any path from $(o,0)$ is
$\sum_n\lambda^n T_o^n$. If $\lambda M_w<1$ then this sum is
finite, hence the number of reproduction trails is a.s.~finite and there
is no weak survival.
\end{proof}

\begin{cor}\label{cor:1}
\begin{enumerate}[1.]
\item For every multigraph, if $M_w=M_s$ there is no pure weak survival.
\item Let $(X,E(X))$ be a non-oriented multigraph.
If $|B(x,n)|^{1/n} \to 1$ for some (equivalently for all) $x
\in X$ then there is no pure weak survival.
\end{enumerate}
\end{cor}
\begin{proof}
\begin{enumerate}[1.]
\item It follows from $1/M_w \le \lambda_w \le \lambda_s =1/M_s$.
\item It is enough to prove that
$M_w\le M_s$. Note that, by the Cauchy-Schwarz inequality,
\[
M_s^{2n}\ge\gamma^{2n}_{x,x}=\sum_{y \in X}\gamma^n_{x,y}\gamma^n_{y,x}=
\sum_{y \in B(x,n)} (\gamma^n_{x,y})^2\ge
\frac{\left(\sum_y\gamma^n_{x,y}\right)^2}{|B(x,n)|}=
\frac{\left(T_x^n\right)^2}{|B(x,n)|},
\]
hence
\[
M_s\ge\limsup_n\sqrt[2n]{\frac{\left(T_x^n\right)^2}{|B(x,n)|}}=
\limsup_n\sqrt[n]{T_x^n}=M_w.
\]
\end{enumerate}
\end{proof}

Let us consider now the question whether
$\lambda_w=1/M_w$. The following theorem states that this equality holds
if the multigraph satisfies a geometrical condition.
By definition of $M_w$, for all fixed $\varepsilon>0$ and $x \in X$, there exists
$n_x$ such that $\sqrt[n_x]{T_x^{n_x}}\ge M_w-\varepsilon$.
We say that
$M_w$ is \textit{attained uniformly} if for all $\varepsilon>0$ there
exists $\bar n=\bar n(\varepsilon)$  for which, for all $x\in X$, $\sup_{n\le \bar n}
\sqrt[n]{T_x^{n}}\ge M_w-\varepsilon$.

\begin{teo}\label{th:conditioncones}
If $(X,E(X))$ is a multigraph such that $M_w$ is attained uniformly
then $\lambda_w=1/M_w$.
\end{teo}

\begin{proof}
Fix $\varepsilon >0$ and $\lambda$ such that
$\lambda (M_w-\varepsilon)>1$. We associate to the BRW
a generalized branching process where the type of each particle
is the site where it is born (although in Section~\ref{sec:gbp} the type
was indexed by $\N$ this is not a restriction since $X$ is at most countable). 
For all $x\in X$ define $n_x$ to be the smallest positive integer such that
$\sqrt[n_x]{T_x^{n_x}} \geq M_w-\varepsilon$.
In this generalized branching process the ``children'' of the initial 
particle (which represents the root of the tree of the process) are
all the particles associated with trails of length $n_o$ starting from $o$. 
Each of these trails ends on a specific vertex in
$B(o,n_o)$, which represents the type of the children generated
there. The offspring distribution $\mu_o$ 
 is supported on $\cE_o:=\{f \in \cE: f(x)=0, \, \forall x \not \in B(o,n_o)\}$
and satisfies
\[
\mu_o(f)=
\pr
\left (
\bigcap_{y \in B(o,n_o)} \{ K_y=f(y)
\}
\right ), \qquad \forall f \in \cE_o,
\]
where $\pr$ is the probability on the space where the BRW is defined
and $K_y$ is the (random) number of trails
of length $n_o$ starting at $o$ and ending at $y$.
The corresponding $\nu_o$ is supercritical in the sense that
\[
G^\prime_o (1) \equiv \sum_{n=0}^\infty n \nu_o(n) = \lambda^{n_o}{T_o^{n_o}}> 1.
\]
This means that $G_o$ has a fixed point $\delta_o<1$.

Analogously, we repeat this
construction for any particle at any site $x$.
The children of such a particle $p$ are the the particles associated with
trails which are prolongations of the trail of $p$ and the difference
between the lengths of the prolongation and of the trail of $p$ is $n_x$.
Clearly  the offspring distribution $\mu_x$ is supported on 
$\cE_x:=\{f \in \cE: f(z)=0, \, \forall z \not \in B(x,n_x)\}$ and is defined as
\[
\mu_x(f)=
\pr
\left (
\bigcap_{y \in B(x,n_x)} \{ K_y=f(y)
\}
\right ), \qquad \forall f \in \cE_x,
\]
where $K_y$ is the (random) number of prolongations, ending at $y$,
of the trail of $p$, such that the difference
between the lengths of the prolongation and of the trail of $p$ is $n_x$.
By Markov property, these laws do not depend on the particle, but only on the
site $x$, hence the definition is well posed.
More precisely, 
$\mu_x$ depends only on the submultigraph $B(x,n_x)$. We call $G_x$
the generating function of $\nu_x$.

These generating
functions $G_x$ are taken from a finite set of $G$'s; indeed in a
bounded degree multigraph the set of the equivalence classes up
to isometries of the balls of radius at most $\bar n$ is finite.
Since all these generating functions are convex, we may apply
Lemma~\ref{th:powerhouse} with $\delta=\max \{\delta_x:x \in X\}$ obtaining
that the generalized branching
process is supercritical.
Since for each $x \in X$ we consider only the particles generated along a path of length $n_x$ (starting
from $x$)
the generalized branching process is dominated by the total number of particles of the original
BRW, hence this last one is supercritical as well. Since $\varepsilon$ was arbitrary,
we deduce that $\lambda_w \le 1/M_w$. Lemma \ref{lem:general} yields the conclusion.
\end{proof}

A large family of multigraphs for which the former condition holds is described by the following theorem.

\begin{teo}\label{th:Tnk}
Let $(X,E(X))$ be a multigraph; let us suppose that there exists $x_0 \in X$, $Y \subseteq X$ and $n_0 \in \N$
such that
\begin{enumerate}[(1)]
\item for all $x \in X$ we have that $B(x,n_0) \cap Y\not = \emptyset$;
\item for all $y \in Y$ there exists an injective map $\varphi_y:X \to X$, such that $\varphi_y(x_0)=y$ and
$n_{\varphi_y(x)\varphi_y(z)} \geq n_{xz}$ for all $x,z \in X$.
\end{enumerate}
Then $M_w$ is attained uniformly and $\lambda_w=1/M_w$.
\end{teo}

\begin{proof}
We fix $\varepsilon >0$. For any given $x  \in X$, condition (1)  implies the existence of $y \in Y$ such that
$\rho(x,y) \leq n_0$, hence $T_y^n\leq T_x^{n+n_0}$ for all $n \in \N$. Using condition (2), we have that
$T_y^n\geq T_{x_0}^n$ for all $n \in \N$, which in turn implies $T_x^{n+n_0} \geq T_{x_0}^n$.
Since $\limsup_{n \to \infty} (T_{x_0}^n)^{1/(n+n_0)}=M_w$, we may find $n_1 \in \N$ such that
$(T_{x_0}^{n_1})^{1/(n_1+n_0)} \ge M_w-\varepsilon$, whence $\overline n(\varepsilon):=n_1+n_0$ satisfies
the hypotheses of Theorem \ref{th:conditioncones}.
\end{proof}

\noindent For a nontrivial example of trees satisfying the hypotheses of the previous theorem
see Example~\ref{exmp0}.
Another important class of multigraphs where $M_w$ is attained uniformly is
described by the following definition (see also Theorem \ref{th:almostregular}).

\begin{defn}\label{Fclass}
Let $(X,E(X))$ and $(Y,E(Y))$ be two multigraphs. A map $\varphi:X
\to Y$ is called a local isomorphism from $X$ onto $Y$ if and only
if
\begin{enumerate}
\item it is surjective,

\item
for all $x \in X$, $y \in Y$
we have
$\sum_{z \in X : \varphi(z)=y}n^X_{xz}=n^Y_{\varphi(x)y}$.
\end{enumerate}
We say that a multigraph (resp.~a graph) $(X,E(X))$ is an $\cF$-multigraph
(resp.~an $\cF$-graph) if
it is locally isomorphic to a finite multigraph (resp.~a finite graph) .
\end{defn}

%

%

Note that a local isomorphism from $X$ to $Y$ does not implies the existence of
a local isomorphism from $Y$ to $X$.
Moreover it is easy to show that, for any local isomorphism, $\varphi(D_X(x)) = D_Y(\varphi (x))$
%
%
%
and that
\begin{equation}\label{eq:XY}
\begin{split}
\sum_{z \in X: \varphi(z)=y}
\gamma^n_{x,z}&=\widetilde \gamma^n_{\varphi(x),y},
\qquad \forall x \in X, \forall y \in Y, \\
T_x^n(X)&=T_{\varphi(x)}^n(Y), \qquad \forall x \in X,
\end{split}
\end{equation}
where $\widetilde \gamma$ refers to paths in $Y$.
The second equation in\eqref{eq:XY} is implied by the first one, which may be proved by induction
using the properties of
$\varphi$.
We note that both quasi-transitive graphs and regular
graphs are $\cF$-multigraphs.
Indeed if $X$ is a quasi-transitive graph, one takes $Y$ as the quotient space
with respect to the action of the automorphism group, $\varphi$ as the quotient map
and $n_{yy^\prime}:=\left|\varphi^{-1}(y^\prime\cap D_X(x)\right|$ where $\varphi(x)=y$
(this definition does not depend on the choice of $x$);
regular graphs of degree $k$ may be mapped on the one-point
multigraph with $k$ loops (and $M_w=k$).
Nevertheless this class contains graphs which are neither regular
nor quasi transitive (see Examples \ref{exmp1}, \ref{exmp1bis} and \ref{exmp0});
moreover the ``regularity'' of $\cF$-multigraphs is only ``local'', indeed one
can easily construct examples of quite irregular $\cF$-graphs.

The following lemma gives a sufficient condition for a graph
to be an $\cF$-graph.

\begin{lem}\label{lem:deg}
Let us consider a graph $(X,E(X))$ such that for all $x,y \in X$ with
$\deg(x)=\deg(y)$ we have that
\[
|\{z \in D(x): \deg(z)=j\}| =|\{z \in D(y): \deg(z)=j\}|, \qquad
\forall j=1,\ldots,M.
\]
Then $(X,E(X))$ is an $\cF$-graph.
\end{lem}
\begin{proof}
Take $Y:=\{i \in \N: \exists x \in X, \deg(x)=i\}$, $n_{ij}:=|\{z \in
D(x): \deg(z)=j\}|$ for some $x \in X$ such that $\deg(x)=i$ (the definition does not
depend on $x$) and $\varphi:=\deg$.
\end{proof}

The following proposition shows how $M_w$, $\lambda_w$ and $\lambda_s$ (or equivalently $M_s$)
are affected by the action of a local isomorphism.

\begin{pro}\label{th:twoprocesses}
Let $(X,E(X))$ and $(Y,E(Y))$ be two connected multigraphs and
suppose that there exists a
local isomorphism $\varphi$ from $X$ onto $Y$. The following assertions hold
\begin{enumerate}
\item $\lambda_w(X)=\lambda_w(Y)$.
\item $\lambda_s(X)\ge\lambda_s(Y)$. If there exists $y \in Y$
such that $|\varphi^{-1}(y)| <+\infty$ then
$\lambda_s(X)=\lambda_s(Y)$.
\item $M_w(X)=M_w(Y)$.
\item $M_w(Y)$ is attained uniformly if and only if
$M_w(X)$ is attained uniformly.
\end{enumerate}
\end{pro}

\begin{proof}
\begin{enumerate}
\item Let $\eta_t$ be a $\lambda$-BRW process on $X$
starting with one particle at site $x$. One may easily show
that
\[
\xi_t(y):= \sum_{x \in \varphi^{-1}(y)} \eta_t(x)
\]
is a $\lambda$-BRW process on $Y$ starting with one
particle at site $\varphi(x)$. It is clear that $\eta_t$ survives
globally if and only if $\xi_t$ does; this implies
$\lambda_w(X)=\lambda_w(Y)$.
\item If $\eta_t$ survives locally
then $\xi_t$ does; hence $\lambda_s(X)\ge\lambda_s(Y)$.
On the other hand, given that $|\varphi^{-1}(y)| <+\infty$, if we start the process $\eta_t$ with one
particle at a site $x \in \varphi^{-1}(y)$ and $\xi_t$ survives
locally (in $y$) the same must be true for $\eta_t$ at some $z \in
\varphi^{-1}(y)$ and hence at $x \in X$.
\item This is a simple consequence of the
equality $T_x^n(X)=T_{\varphi(x)}^n(Y)$ which holds for all $x \in
X$ and $n \in \N$.
\item It follows from the facts that $M_w(X)=M_w(Y)$ and $T_x^n(X)=T_{\varphi(x)}^n(Y)$.
\end{enumerate}
\end{proof}

We note that, according to the previous proposition,
if $(X,E(X))$ is locally isomorphic to a multigraph $(Y,E(Y))$ which satisfies the
hypotheses of Theorem \ref{th:Tnk}, then the same conclusions of this theorem hold for $(X,E(X))$.
In particular if $(Y,E(Y))$ is a finite multigraph then $M_w(Y)=M_s(Y)$ and
$\lambda_w(Y)=\lambda_s(Y)$.

\begin{teo}\label{th:almostregular}
Let $(X,E(X))$ be an $\cF$-multigraph, then $M_w$ is attained uniformly and $\lambda_w=1/M_w$. Moreover
if $\lambda=\lambda_w$ the $\lambda$-BRW on $X$ dies out globally almost surely.
\end{teo}
\begin{proof}
Let $(X,E(X))$ be locally isomorphic to the finite multigraph $(Y,E(Y))$.
We note that $M_w(Y)$ is attained uniformly (since $Y$ is finite) whence,
by Proposition \ref{th:twoprocesses}, $M_w(X)$ is attained uniformly.

Since the global behavior of the $\lambda$-BRW $\eta_t$ on $X$ is the same
as the corresponding behavior of the induced $\lambda$-BRW $\xi_t$ on $Y$
(see the proof of Proposition
\ref{th:twoprocesses}), then
Theorem \ref{th:conditioncones} and Proposition \ref{th:twoprocesses} imply $\lambda_w(X)=1/M_w(X)=\lambda_w(Y)=\lambda_s(Y)$. By
Theorem \ref{th:pemantleimproved} each $\lambda_s$-BRW dies out
locally a.s.; moreover, since $Y$ is a finite multigraph,
$\xi_t$  dies out globally a.s., hence the same holds for $\eta_t$.
\end{proof}

\begin{rem}\label{rem:finitegraph}
{\rm
It is natural to wonder how $M_s$, $M_w$, $\lambda_s$ and $\lambda_w$ are affected
by local modifications of the multigraphs $(X,E(X))$ (such as, for instance,
attaching a complete finite, graph to a vertex of $X$ or removing a set of
vertices and/or edges).

If $(X,E(X))$, $(Y,E(Y))$ are two multigraphs and $\psi:Y \to X$
is an injective map such that $n_{\psi(x)\psi(y)} \ge \widetilde n_{xy}$ for
all $x,y \in Y$ (where $\widetilde n$ refers to $Y$) then
$\lambda_w(X) \le \lambda_w(Y)$,
$\lambda_s(X) \le \lambda_s(Y)$,
$M_w(X) \ge M_w(Y)$, $M_s(X) \ge M_s(Y)$.

In certain cases it is easy to show that  the existence of a pure weak phase on
$X$ implies the existence of a pure weak phase on some submultigraph; indeed
if $Y$ is a finite subset of $X$ such that $X \setminus Y$ is divided into
a finite number of connected multigraphs $X_1, \ldots, X_n$ (which is certainly true if
$n_{xy}>0$ is equivalent to $n_{yx} >0$ for all $x,y \in X \setminus Y$), then for every
$\lambda \in (\lambda_w(X),\lambda_s(X))$ the $\lambda$-BRW leaves eventually a.s.~the subset $Y$.
Hence it survives (globally but not locally) at least on one connected component; this means that,
although
$\lambda_s(X_i) \ge \lambda_s(X)$, $\lambda_w(X_i) \ge \lambda_w(X)$ for all $i=1, \ldots,n$,
there exists $i_0$ such that $\lambda_w(X_{i_0}) = \lambda_w(X)$. The existence of a pure
weak phase on $X_{i_0}$ follows from
$\lambda_s(X_{i_0}) \geq \lambda_s(X)>\lambda_w(X)=\lambda_w(X_{i_0})$.

Moreover if there exists a subset $Y$ as above such that $\lambda_w(X_i) > \lambda_w(X)$ for all $i$, 
then there is no pure weak phase for the BRW on $X$.
Take for instance
a graph $(X^\prime,E(X^\prime))$ and $k \in \N$ such that
$1/k < \lambda_w(X^\prime)$. Attach a complete graph of degree $k$ to a vertex of
$X^\prime$, we obtain a new graph $X$ such that $\lambda_s(X)=\lambda_w(X) 
\leq 1/k < \lambda_w(X^\prime)$; hence even if the BRW on $X^\prime$ has a pure weak phase,
the BRW on $X$ has none.}
\end{rem}

\subsection{Nonamenability and weak phase}\label{sec:nonamenability}

In this section we consider only non oriented multigraphs.
A multigraph
$(X,E(X))$ is \textit{nonamenable} if
\[
\inf
\left \{
\frac{|\partial_E(S)|}{|S|}: S \subseteq X, |S| < \infty
\right \}=:\iota_X>0,
\]
where $\partial_E(S)$ is the set of edges $(x,y,i)\in E(X)$ such that $x\in
S$ and $y\not\in S$.

We define $N:l^2(X) \to l^2(X)$ by $Nf(x):=\sum_{y \in X} n_{xy}
f(y)$ which is a bounded, linear operator with
$\|N\| \leq M$.
It is well known that on a regular, non-oriented graph (where $M=M_w$) the
existence of the weak phase is equivalent to nonamenability (see \cite[Theorem 2.4]{cf:Stacey03}).
Indeed on regular,  non-oriented graphs $M_s<M_w$ is equivalent to nonamenability:
one easily proves that $M_s=\|N\|$ (see Lemma \ref{lem:norm} and the reference therein); moreover
$\|N\|=M \|P\|$ where $P$ is the transition operator associated to the simple random walk and
Gerl proved that $\|P\| < 1$ is equivalent to nonamenability (see \cite{cf:Gerl}).
Hence using 
Theorem~\ref{th:pemantleimproved} and Theorem~\ref{th:almostregular} we obtain an alternative proof of
\cite[Theorem 2.4]{cf:Stacey03}. 

Now we show that, for non-oriented $\cF$-multigraphs,
nonamenability is equivalent to the existence of a pure weak phase which, in turn,
is equivalent to $M_s<M_w$.

\begin{lem}\label{lem:norm}
$N$ is self
adjoint and $\|N\|=\rho(N)=M_s$ where $\rho(N)=\lim_{n \to \infty}
\|N^n\|^{1/n}$ is the spectral radius of $N$.
\end{lem}
\begin{proof}
The self-adjointness of $N$ is easy and $\|N\|=\rho(N)$ is a
standard property which follows from the Spectral Theorem for any
normal (hence self-adjoint) operator. To prove that $M_s=\|N\|$ one
proceeds essentially as in \cite[Lemma 2.2]{cf:Stacey03}.
\end{proof}

\noindent
The following theorem implies the analogous results for regular
and quasi-transitive graphs.

\begin{teo}\label{th:nonam}
Let $(X,E(X))$ be a non-oriented $\cF$-multigraph.
Then $\lambda_w<\lambda_s$ if
and only if $(X,E(X))$ is nonamenable.
\end{teo}
Before proving this statement, we need a technical result concerning
the Dirichlet norm of $l^2$ functions.
Given $f\in l^2(X)$,
define
\[
\Vert f\Vert_{D(2)} = \left( \sum_{x,y\in X} n_{xy}|f(x)-f(y)|^2
\right)^{1/2}.
\]

\begin{lem}\label{th:norms}
Let $(X,E(X))$ be a nonamenable multigraph. Then
there exists $c>0$ such that, for all $f\in l^2(X)$,
\[
\Vert f\Vert_{D(2)}\ge c\Vert f\Vert_2.
\]
\end{lem}
\begin{proof}
The proof is analogous to the one of
\cite[Theorem 2.6]{cf:Stacey03} (one has to deal carefully
with the presence of $n_{xy}$), hence we omit it.
\end{proof}

\begin{proof}[Proof of Theorem \ref{th:nonam}]
We follow the proof of \cite[Theorem 3.1]{cf:Stacey03}.
Let $(X;E(X))$ be nonamenable, $N=(n_{xy})_{x,y\in X}$ be its adjacency
matrix and $\widetilde N=(\tilde n_{xy})_{x,y\in Y}$ be the
adjacency matrix of the finite multigraph $(Y,E(Y))$ which $(X,E(X))$ is locally isomorphic to.
We must prove that $\Vert
N\Vert<M_w$. By definition of local isomorphism we have that
$\tilde n_{\varphi(x)\varphi(y)}= \sum_{z:\varphi(z)=\varphi(y)} n_{xz}$. By the Perron-Frobenius
theorem $\widetilde N$ has largest positive eigenvalue $M_w$ with
associated positive eigenvector $(a_1,\ldots, a_k)$ ($k$ being the
cardinality of $Y$). Then
\begin{equation}
\label{eq:eigenv} M_wa_{\varphi(x)}=\sum_{y^\prime\in
Y}\widetilde n_{\varphi(x)y^\prime}a_{y^\prime}=\sum_{y^\prime\in
Y}\sum_{y\in\varphi^{-1}(y^\prime)} n_{xy}a_{\varphi(y)}= \sum_{y\in X}n_{xy}a_{\varphi(y)}.
\end{equation}
Take $f\in l^2(X)$. Applying equation~\eqref{eq:eigenv} and the
fact that $(X,E(X))$ is non-oriented
\[
\begin{split}
M_w^2\Vert f\Vert_2^2 & = M_w^2\sum_{y\in X} (f(y))^2\\
& = M_w\sum_{y\in X} \left(\sum_{x \in X}
n_{yx}\frac{a_{\varphi(x)}}{a_{\varphi(y)}}
\right)(f(y))^2\\
& = M_w\sum_{x\in X} \left(\sum_{y \in X}
n_{yx}\frac{a_{\varphi(x)}}{a_{\varphi(y)}}
\right)(f(y))^2\\
& = \sum_{x\in X} M_w a_{\varphi(x)}\sum_{y \in X}
n_{yx}\frac{(f(y))^2}{a_{\varphi(y)}}
\\
& = \sum_{x\in X}\left(\sum_{z \in X} n_{xz}a_{\varphi(z)}\right)
\left(\sum_{y \in X} n_{xy}\frac{(f(y))^2}{a_{\varphi(y)}}\right).
\end{split}
\]
Hence
\[
\begin{split}
M_w^2\Vert f\Vert_2^2 -\Vert Nf\Vert_2^2 & = \sum_{x\in X}
\sum_{z,y \in X} n_{xz}n_{xy}\left[
\frac{a_{\varphi(z)}}{a_{\varphi(y)}}(f(y))^2- f(z)f(y)\right]
\\
& = \frac12 \sum_{x\in X} \sum_{z,y \in X}
n_{xz}n_{xy}a_{\varphi(y)}a_{\varphi(z)} \left[
\frac{f(y)}{a_{\varphi(y)}}- \frac{f(z)}{a_{\varphi(z)}}
\right]^2\\
&\ge \frac12 (\min a_i)^2 \Vert g \Vert_{D(2)},
\end{split}
\]
where $g(x)=f(x)/a_{\varphi(x)}$ is considered as a map on the
multigraph $G^2=(X,\bar E(X))$ with adjacency matrix $\bar N$ defined by
$\bar n_{xy}=\sum_{z \in X} n_{zx}n_{zy}$ and
$\bar E(X):=\{(y,z,i) : 1 \le i \le \bar n_{xy}\}$. Applying
Lemma~\ref{th:norms} to each connected component of $G^2$ (note that each of them is nonamenable) and noting that $\Vert g\Vert_2^2\ge
D^2\Vert f\Vert_2^2$ for $D^{-1}=\max(a_i)$, we have that for some
$C>0$
$$
M_w^2\Vert f\Vert_2^2-\Vert Nf\Vert_2^2\ge C\Vert f\Vert_2^2,
$$
whence $\Vert N\Vert\le\sqrt{M_w^2-C}<M_w$.
\smallskip

Suppose now that $(X,E(X))$ is amenable and fix $\eps>0$. Then for some
finite set $S\subset X$, $|\partial_E S|/|S|<\eps$. Define
$f(x)=a_{\varphi(x)}\ident_{S}(x)$. If $x\in S$ and
$D(x)\cap S^c=\emptyset$, then, by \eqref{eq:eigenv}, $Nf(x)=M_wf(x)$. Hence
$$
\Vert Nf\Vert_2^2\ge M_w^2\left( \Vert f\Vert_2^2-\eps |S|(\max
a_i)^2 \right),
$$
and
$$
\frac{\Vert Nf\Vert_2^2}{\Vert f\Vert_2^2}\ge M_w^2-
\eps\left(\frac{\max a_i}{\min a_i}\right)^2.
$$
By taking $\eps$ arbitrarily small we prove that $\Vert
N\Vert \ge M_w$, whence $M_s=M_w$ (recall that $\|N\|=M_s \le M_w$).
\end{proof}


\subsection{Examples}\label{sec:examples}

The first two explicit examples listed hereafter show that
the class of $\cF$-multigraphs is larger than the union of regular
and quasi-transitive multigraphs.
Both these examples are modifications of regular graphs:
Example \ref{exmp1} is obtained by attaching an edge to each vertex,
Example \ref{exmp1bis} by drawing a ``bridge with intermediate station'' between some of the vertices.

\begin{exm}\label{exmp1}
{\rm \noindent
Take a square and attach to every vertex a branch of a homogeneous tree of degree $3$, obtaining a regular graph (of degree $3$)
which is not quasi transitive. If we attach now to each vertex a new edge with a new endpoint we obtain a
non-oriented, nonamenable
$\cF$-graph $(X,E(X))$ which is neither regular nor quasi transitive. It is easily seen
(by Lemma \ref{lem:deg}) to be locally isomorphic to a multigraph with
adjacency matrix
\[
N=
\begin{pmatrix}
3 & 1\\
1 & 0 \\
\end{pmatrix}.
\]
According to Theorem \ref{th:nonam}, the BRW on this graph has a pure weak phase.
}
\end{exm}

\begin{exm}\label{exmp1bis}
{\rm \noindent
Take an infinite graph $(X,E(X))$ with set of vertices $X=\{x_1,x_2,\ldots\}$.
If $Y=\{y_1,y_2,\ldots\}$ is another countable set, disjoint from $X$, we may consider the graph with set of vertices
$Z:=X \cup Y$ and
\[
E(Z):=E(X)\cup \bigcup_{i=1}^\infty \{(x_{2i-1},y_i), (y_i,x_{2i-1}), (x_{2i},y_i), (y_i,x_{2i})\};
\]
roughly speaking we join $x_{2i-1}$ and $x_{2i}$ by a bridge and we cut this bridge into two edges by using a new vertex
$y_i$.
If the graph $X$ is nonamenable then it is possible to show that the (multi)graph $Z$ is nonamenable as well. By choosing
$(X,E(X))$ regular (with $\deg \equiv k$) we obtain an $\cF$-graph which
(by Lemma \ref{lem:deg}) is locally isomorphic to a multigraph with adjacency matrix
\[
N=
\begin{pmatrix}
k & 1\\
2 & 0 \\
\end{pmatrix}.
\]
Again, by choosing accurately $(X,E(X))$
and ordering wisely its vertices we may obtain a graph which is neither quasi transitive nor regular.
}
\end{exm}
The following trees are natural examples of graphs which are not quasi transitive and, nevertheless, are not
``too irregular''. We show that, for these trees, nonamenability is equivalent to the existence of a pure
weak phase and the proof is not a direct application of
Theorem \ref{th:nonam}.

\begin{exm}\label{exmp0}
{\rm
Given a sequence of positive natural numbers $\{n_k\}_k$
we construct a non-oriented, rooted tree $\T$ (with root $o$) such that if $x\in \T$ satisfies $\rho(o,x)=k$
then $\deg(x)=n_k+1$. We call this radial graph $T_{\{n_k\}}$-tree.
If the sequence is periodical of period $d$, then Theorem \ref{th:Tnk} applies
with $x_0=o$, $n_0=d$, $Y:=\cup_{n\in\N} B(o,nd)$ and  $\varphi_y$
(where $y \in Y$) maps isomorphically
the tree $\T$ onto the subtree branching from $y$. We call $T_i$
the $T_{\{n^\prime_k\}}$-tree obtained by means of this construction where $n^\prime_k:=n_{k+i}$.
Roughly speaking we construct $T_1, \ldots, T_{d-1}$,
by using cyclic permutations of the sequence $\{n_1, \ldots, n_{d-1}\}$. Obviously $\T=T_1$.
Since $T_i$ may be mapped into $T_j$ for all $i,j=1, \ldots, d$ (in the sense of Remark \ref{rem:finitegraph})
then $\lambda_w(T_i)$, $\lambda_s(T_i)$, $M_w(T_i)$, and $M_s(T_i)$ do not depend on $i$.

Let us consider the finite cyclic graph $\widetilde Y:=\{y_1, \ldots, y_d\}$ where $n_{y_i y_{i+1}}=n_{y_{i+1} y_i}=1$
for all $i=1,\ldots d$ (with the identification $y_{d+1} \equiv y_1$). To each vertex $y_i$ we attach
$n_i-1$ copies of $T_{i+1}$ (again with the identification $T_{d+1}=T_1$), each of them by using a two-way edge.
We denote this  connected, non-oriented $\cF$-graph  by $(X,E(X))$; indeed it may be mapped
onto the finite multigraph $Y^\prime$ where $Y^\prime=\widetilde Y$ and $n^\prime_{y_i y_{i+1}}=n_i$, $n^\prime_{y_{i+1} y_i}=1$
for all $i=1,\ldots d$. Note that $X$ is neither quasi transitive nor regular, unless $n_i = 1$ for all $i$.
$X$ is nonamenable if and only if $T_1$ is nonamenable, that is, if and only if there exists $i$ such that
$n_i\ge 2$. In this case, according to Theorem \ref{th:nonam}, $\lambda_w(X)< \lambda_s(X)$, hence by Remark~\ref{rem:finitegraph}
(considering $X \setminus \widetilde Y$) there exists
$i$ such that $\lambda_w(T_i)< \lambda_s(T_i)$.
This means that for all $i$ we have $\lambda_w(T_i)<  \lambda_s(T_i)$ and there is a pure
weak phase on $T_i$. On the other hand, if $n_i\equiv 1$ for all $i=1,\ldots,d$, then there is no pure weak phase (Corollary \ref{cor:1}).
}
\end{exm}

\section{Modified BRW and BRW on weighted graphs} \label{sec:modif}

\subsection{Modified BRW} \label{sec:MBRW}

In this section we consider an irreducible random walk $(X,P)$. In
the case of simple random walks some
of the results of this section may be found also in \cite{cf:Stacey03}.
We
study the modified BRW where each particle at
site $x$ dies at rate 1 and breeds at rate $\lambda$ and sends the
offspring randomly according to the probability distribution
$p(x,\cdot)$.

We denote by $p^{(n)}(x,y)$ the $n$-step transition probabilities
from $x$ to $y$ ($n\ge0$) and by $f^{(n)}(x,y)$ the probability
that the random walk starting from $x$ hits $y$ for the first time
after $n$ steps ($n\ge 1$). Then we define the corresponding
generating functions $G(x,x|z)=\sum_{n\ge0}p^{(n)}(x,x) z^n$
and $F(x,x|z)=\sum_{n\ge1}f^{(n)}(x,x) z^n$, where
$x\in X$, $z\in\C$ (further details can be found in
\cite[Chapter I.1.B]{Woess2}, where $F$ is called $U$).

The expected number of trails along a path
$\Pi=\{x_0,\ldots,x_n\}$  is equal to
$\lambda^n \prod_{i=0}^{n-1} p(x_i,x_{i+1})$. Hence the expected
number of trails  along paths starting from
$x$ and reaching $y$ for the first time is equal to
$F(x,y|\lambda)$. If $x$ is equal to $y$ we call them first
generation trails in $x$.
Since
$G(x,x|\lambda)={1}/{(1-F(x,x|\lambda))}$
and the radius $R$ of $G$ does not depend on the choice of $x$, we
have that $R=\max\{\lambda:F(x,x|\lambda)\le 1\}$.

\begin{teo}\label{th:modif}
For the modified BRW $\lambda_w=1$ and if $\lambda=1$ there is
global extinction almost surely. Moreover $\lambda_s=R$ and if $\lambda=R$
there is local extinction almost surely.
\end{teo}
\begin{proof}
The total number of particles $T_t$ is a branching process with
rate $\lambda$, whence the claim for $\lambda_w$ follows.
As for the second claim, the proof is the same as in Theorem
\ref{th:pemantleimproved} using $F$ instead of $\Phi$.
\end{proof}

\noindent The following Corollary is the analog of Theorem \ref{th:nonam} (see \cite{cf:Gerl}
for the definition of \textit{strongly reversible random walk}).

\begin{cor}\label{cor:modif}
For the modified BRW, the existence of a pure weak phase is equivalent to $R>1$. If $P$ is a strongly reversible random walk
then the existence of the pure weak phase is equivalent to nonamenability.
\end{cor}
\begin{proof}
The result is a simple consequence of Theorem \ref{th:modif} and the main theorem of \cite{cf:Gerl}.
\end{proof}

\subsection{BRW on weighted graphs}

Our methods apply, with minor modifications, to more general BRWs, which
generalize simultaneously BRWs on multigraphs and modified BRWs.

Let us consider $(X,N)$ where $X$ is a countable (or finite) set and $N=(n_{xy})_{x,y \in X}$ is
a matrix of nonnegative \textit{weights} (that is, $n_{xy} \ge 0$) such that $\sup_{x \in X} \sum_{y \in X}
n_{xy} = M< \infty$. We suppose that $N$ is irreducible in the sense that
$(X,E(X))$, where $E(X):=\{(x,y) \in X \times X: n_{xy}>0 \}$, is a connected graph. We call $(X,N)$ a \textit{weighted graph}.

The $\lambda$-BRW is defined by setting the reproduction rate on every edge
$(x,y)$ as $\lambda n_{xy}$; hence, to each path $\{x_0, \ldots, x_n\}$ there corresponds a weight $\prod_{i=0}^{n-1} n_{x_i x_{i+1}}$. We
define $\gamma^n_{x,y}$, $T^n_x$, $\phi^n_{x,y}$, $M_s$ and $M_w$ as in Section \ref{subse:def}.

It is clear that the BRW on multigraphs and the modified (according to an irreducible random
walk) BRW may be viewed as BRWs on weighted graphs.
Moreover the expected number of trails along a path
$\{x_0, x_1, \ldots, x_n\}$ is
$\lambda^n \prod_{i=0}^{n-1} n_{x_i x_{i+1}}$.
Substituting the word ``multigraphs'' with ``weighted graphs'' all the results of Sections \ref{sec:def}, \ref{sec:main} and
\ref{sec:MBRW} still hold (with the exception of Theorems \ref{th:conditioncones} and
\ref{th:Tnk}) with unimportant modifications.
In particular, extending Definition \ref{Fclass} verbatim to weighted graphs, one can prove Theorem \ref{th:almostregular}, since in this
case $n_{xy}$ may take just a finite number of values and it is possible to apply Lemma \ref{th:powerhouse} as we did in
Theorem \ref{th:conditioncones}.
For \textit{regular} weighted graphs (that is, $\sum_{y \in X}
n_{xy} = M$ for all $x \in X$) one proves results analogous to the ones of Section \ref{sec:MBRW}.

%
%
%
%
%
%

\section{Open questions}
\label{sec:fr}

As we stated in Section \ref{sec:intro}, this paper is motivated
by three main issues: the identification of the critical value $\lambda_w$, the behavior of the process
when $\lambda=\lambda_s$ or $\lambda=\lambda_w$ and the existence of the pure weak phase.

To complete the first point one should verify whether the equality $\lambda_w=1/M_w$ holds for every
multigraph or if $M_w$ characterizes the critical value $\lambda_w$ only on a restricted class of multigraphs.

As for the second one, the open question is the following: is it possible to construct a multigraph where if $\lambda=\lambda_w$
the process does not die out globally? In particular, is it possible to find
a multigraph where $\lambda_s=\lambda_w$ but the $\lambda_w$-BRW does not die out
globally (it certainly does locally)?

Finally, dealing with the existence of a pure weak phase, it is well known that there is
no equivalence, in general, with nonamenability. We proved that this equivalence holds, for instance, for the class
of non-oriented $\cF$-multigraphs; we do not know what can be said in the case of oriented $\cF$-multigraphs.
To be precise: is there a nonamenable, oriented $\cF$-multigraph,
where the BRW has no weak phase?
On the other hand, is it possible to find an amenable, oriented $\cF$-multigraph where $\lambda_s=\lambda_w$?

\end{document}